\documentclass{article}
\title{A prime prime primer}
\author{Michael Caola\\6 Normanton Rd., Bristol, BS8 2TY, U.K.\\caola@blueyonder.co.uk}
\begin{document}

\maketitle
\begin{abstract}
We present a simple, closed formula which gives all the primes in order. It is a simple product of integer floor and ceiling functions.
\end{abstract}
\section{Introduction}

Interest in prime numbers continues apace, e.g. computing, secure communications, quantum computing, statistical physics, engineering and, of course, mathematics. However, we are still far from finding a useful function $p(n)$ which generates all the primes in order. We present such a function next and discuss its usefulness [1]. All variables below, save $x$, are integers.

\section{Analysis}
Consider the product function

\begin{equation}
p(n)=\prod_{n'=2}^{\lfloor\sqrt{n}\rfloor}\biggl( \bigg\lceil \frac{n}{n'} \bigg\rceil- \bigg\lfloor \frac{n}{n'} \bigg\rfloor \biggr),
\end{equation}
where the integer {\it floor} $\lfloor x \rfloor$  ({\it ceiling} $\lceil x \rceil$) functions  give the greatest (smallest) integer not greater (not smaller) than $x$: e.g.
 $\lfloor 7 \rfloor=\lceil 7 \rceil=7$, $\lfloor 3.7 \rfloor=3$, $\lceil 3.7 \rceil=4$, $\lceil -7.3 \rceil=-7$, $\lfloor -7.3 \rfloor=-8.$
\par

Then, for $n=2,3,4 ... $,

\begin{equation}
p(n)=
\left\{\begin{array}{cl}
1,& \mbox{if n is prime}\\         
0,& \mbox{if not}
\end{array}\right.\
\end{equation}\\

Also, we have directly from (2) that the number $\pi(n)$ of primes not exceeding $n$ is
\begin{equation}
\pi(n)=\sum_{n'=1}^{n}p(n')
\end{equation}

\section{Discussion}
\begin{itemize}
\item
That function (1) satisfies (2) may be seen by inspection of $\lceil\rceil$ - $\lfloor\rfloor$:
\begin{equation}
\lceil x \rceil-\lfloor x \rfloor =
\left\{\begin{array}{cl}
0,& \mbox{if $x$ is an integer}\\         
1,& \mbox{if not}\\
\end{array}\right.\
\end{equation}\\

Thus $p(n)=1$ only if all the divisors $n'$ of $n$ give non-zero remaninders, which is the definition of primality.

\item
What constitutes a ``useful function'' $p(n)$ is partly subjective [1], so we  describe our $p(n)$ as: simple, compact, transparent, and using only two basic notions ({\it floor} $\lfloor \rfloor$) and ({\it ceiling} $\lceil \rceil$) of integer mathematics.
\item
We have programmed and checked algorithm (1,2,3) on a PC.

\end {itemize}

\section*{Reference}
[1] Hardy, G. H. \& Wright, E. M. (1978) 'An introduction to the theory of numbers', ( 5th edition). Clarendon Press, Oxford.

\end{document}